\documentclass{article}
\usepackage[T1]{fontenc}
\usepackage{fancyvrb}
\usepackage{comment}
\usepackage{enumerate}
\usepackage{amssymb}
\usepackage{latexsym}
\usepackage{amsmath}
\usepackage{amsfonts}
\usepackage{amsthm}
\newcommand{\rank}{\operatorname{rank}}
\newcommand{\cR}{\mathcal{R}}
\newcommand{\cO}{\mathcal{O}}
\newcommand{\cS}{\mathcal{S}}
\newcommand{\cA}{\mathcal{A}}
\newcommand{\GF}{\operatorname{GF}}

\newcommand{\cK}{\mathcal K}

\newcommand{\Sz}{\operatorname{Sz}}
\newcommand{\GL}{\operatorname{GL}}
\newcommand{\PG}{\operatorname{PG}}
\newcommand{\PGL}{\operatorname{PGL}}

\newcommand{\keywords}[1]{\leftline{{\bfseries Keywords: }#1}}

\newtheorem{prop}{Proposition}
\newtheorem{theorem}[prop]{Theorem}
\newtheorem{definition}[prop]{Definition}
\theoremstyle{definition}
\newtheorem{remark}[prop]{Remark}

\ifx\smallsetminus\undefined
\def\smallsetminus{\setminus}
\fi

\title{An algorithm for constructing some maximal arcs in $\PG(2,q^2)$}
\author{A. Aguglia ${}^*$ \and L. Giuzzi
  \thanks{Research supported by  the Italian
    Ministry MIUR, Strutture geometriche, combinatoria e loro
    applicazioni.}}

\begin{document}
\maketitle
\begin{abstract}
In 1974, J. Thas constructed a  new class of maximal arcs for the
Desarguesian plane of order $q^2$. The construction relied upon the
existence of a regular spread of tangent lines to an ovoid in
$\PG(3,q)$ and, in particular, it does apply to the Suzuki--Tits
ovoid. In this paper, we describe an algorithm for obtaining a
possible representation of such arcs in $\PG(2,q^2)$.
\end{abstract}
\section{Introduction}
In a finite projective plane of order $q$,
a maximal $(k, n)$--arc
$\mathcal K$, where $k\geq 1$ and $2\leq n \leq q+1$, is a non--empty
set of $k$ points  which is met by every line of the plane
in either $0$ or $n$ points. The integer $n$ is called
the \emph{degree} of the arc $\cK$.

Trivial examples of maximal arcs of degree
$q+1$ and $q$ are respectively the set of
all the points of $\PG(2,q)$ and the set of the points of an affine
subplane $AG(2,q)$ of $\PG(2,q)$.

It has been shown in \cite{BBM, BB2}, that non--trivial
maximal arcs do not exist in $\PG(2,q)$ for $q$ odd. On the other
hand, when  $q$ is even,
 several classes of non--trivial
maximal arcs are known. In fact, hyperovals and their duals are
maximal arcs. Apart from them, known constructions of degree $n$
maximal arcs in $PG(2,q)$ are: one construction by Denniston
\cite{D} based on a linear pencil of conics, two constructions of J.
A. Thas \cite{T74,T80}, constructions by Mathon \cite{M},
and by Hamilton and  Mathon \cite{HM}  utilising closed sets of
conics. However, most of the known examples of degree $n$ maximal
arcs (with the notable exception of a class of maximal arcs arising
from the \cite{T74} construction) consist of the union of $n-1$
pairwise disjoint conics, together with their common nucleus $N$. We
shall call these arcs \emph{conical}. Observe that any conical
maximal arc is covered by a completely reducible curve of degree
$2n-1$, whose components are $n-1$ conics and a line through the
point $N$; in \cite{AGK} it is shown that such a curve has
minimum degree. 
In the present paper we
determine, using the computer algebra package GAP \cite{GAP4},
equations for algebraic plane curves
of minimum degree passing through all the points of a maximal arc
$\cK$; in particular, we are interested in those cases in which
$\cK$ is not conical. In particular, we will show that arcs arising
from the Suzuki--Tits ovoids by \cite{T74} construction
cannot be covered by a curve of low degree.

\section{Reguli in  $\PG(3,q)$}
We recall some basic properties of reguli and spreads of $\PG(3,q)$,
see \cite{FPS3d}.
\begin{definition}
A \emph{regulus}
 of $\PG(3,q)$ is a collection of $q+1$ mutually disjoint lines such
 that any line of $PG(3,q)$ meeting three of them necessarily
 meets them all.
\end{definition}
A standard result, see \cite{HT}, shows that any
three pairwise disjoint lines $\ell_1, \ell_2$, $\ell_3$  of
$\PG(3,q)$ lie together in a unique regulus, say
$\cR(\ell_1,\ell_2,\ell_3)$.

\begin{definition}
  Let $\ell_1,\ell_2,\ell_3$ be $3$ pairwise disjoint lines
  of $\PG(3,q)$. The \emph{opposite regulus} to
  $\cR(\ell_1,\ell_2,\ell_3)$
  is the set
  \[ \cR^o(\ell_1,\ell_2,\ell_3) \]
  of all lines $\ell$ of  $\PG(3,q)$ such that
  \[ \ell\cap\ell_i\neq\emptyset, \text{ for $i=1,2,3$.} \]
\end{definition}
The set $\cR^o(\ell_1,\ell_2,\ell_3)$ is also a
regulus.
We may compute
the regulus containing $\ell_1,\ell_2$ and $\ell_3$ as
the set
\[ \cR(\ell_1,\ell_2,\ell_3)=\cR^o(m_1,m_2,m_3), \]
where $m_1,m_2,m_3$ are distinct elements of
$\cR^o(\ell_1,\ell_2,\ell_3)$.

\begin{definition}
  A $k$--span of $\PG(3,q)$ is a set of $k$
  mutually skew lines. A $(q^2+1)$--span
  is called a \emph{spread}.
\end{definition}
Observe that a spread is a partition of the points of
$\PG(3,q)$ in disjoint lines.
\begin{definition}
  A spread $\cS$ is \emph{regular} or
  \emph{Desarguesian}, if for any
  three lines $\ell_1,\ell_2,\ell_3\in\cS$,
  \[ \cR(\ell_1,\ell_2,\ell_3)\subseteq\cS. \]
\end{definition}
Any two regular spreads of
$\PG(3,q)$ are projectively equivalent.
In order to implement the construction of the arc
in Section \ref{sect:3}, we must describe a spread
of tangent lines to a given ovoid.
We shall make use of the following notion of closure.
\begin{definition}
 The \emph{regular closure} of a set $S$ of lines of $\PG(3,q)$
 is the smallest set $T$ of lines of $\PG(3,q)$ containing $S$
 such that for any $3$ distinct elements $\ell_1,\ell_2,\ell_3\in T$,
 \[ \cR(\ell_1,\ell_2,\ell_3)\subseteq T. \]
\end{definition}
Examples of sets closed under this operation are regular spreads of
$\PG(3,q)$ and reguli.
In fact, a regular spread is uniquely determined by four of its
lines, supposed they are in suitable position.
\begin{theorem}
 There exists exactly one regular spread containing
 any given $4$ mutually skew lines $\ell_1,\ell_2,\ell_3,\ell_4$
 of $\PG(3,q)$, provided that $\ell_4\not\in\cR(\ell_1,\ell_2,\ell_3)$.
\end{theorem}
\begin{proof}
By \cite{Po} there is a Desarguesian spread containing any
two reguli with $2$ lines in common.
We now show that this spread is the regular closure
of $\ell_1,\ell_2,\ell_3,\ell_4$.
Any Desarguesian spread containing $\cR(\ell_1,\ell_2,\ell_4)$
and $\cR(\ell_1,\ell_3,\ell_4)$ must clearly
contain also the $(q^2-q+2)$--span of lines given by
\[ \bigcup_{\begin{subarray}{c}
    x\in\cR(\ell_1,\ell_2,\ell_3)\\
    x\neq\ell_1
    \end{subarray}}\cR(\ell_1,x,\ell_4)
\]
By \cite[Lemma 17.6.2]{FPS3d},  a spread containing
such span is unique. The result follows.
\end{proof}

\section{Thas \cite{T74} maximal arcs}
\label{sect:2}
We shall make extensive use of the representation of $\PG(2,q^2)$ in
$\PG(4,q)$ due to Andr\'e \cite{A} and Bruck and Bose
\cite{BBo1,BBo2}.
\par
Let $\PG(4,q)$ be a projective $4$--space over the finite field
$\GF(q)$ and let suppose $\cS$ be a regular spread of a fixed
hyperplane $\Sigma=\PG(3,q)$ of $\PG(4,q)$.
Then $PG(2,q^2)$ can be
represented as the incidence structure whose points are   the points
of $\PG(4,q)\setminus\Sigma$ and the elements of $\cS$, and whose
lines are the planes of $\PG(4,q)\setminus\Sigma$ which meet
$\Sigma$ in a line of $\cS$ and the spread $\cS$. In particular,
$\cS$ represents the ``line at infinity'' of the affine plane
$AG(2,q^2)\subseteq\PG(2,q^2)$.
Recall that 
projectively equivalent spreads of $\PG(3,q)$
induce, via Bruck--Bose construction
isomorphic projective planes or order $q^2$.
In particular, any two regular spreads of $\PG(3,q)$
induce a representation of the Desarguesian projective
plane $\PG(2,q^2)$.
\par
Using the aforementioned  model, Thas obtained
maximal arcs in the Desarguesian plane as follows. Let $\mathcal O $
be an ovoid in the hyperplane $\Sigma$ such that every element of
the spread $\mathcal S$  meets $\cO$ in exactly one point. Fix a
point $V$  in $ PG(4,q)\setminus \Sigma$ and let $\overline{\cK}$ be
the set of points in $ PG(4,q)\setminus \Sigma$ collinear with $V$
and a point on $\cO$. Then $\overline{\cK}$ corresponds to a maximal
$(q^3-q^2+q,q)$--arc $\cK$ in $PG(2,q^2)$.

In \cite{T74}, it has been remarked that if $\cO$ is an elliptic
quadric then the maximal arc thus constructed turns out to be of
Denniston type. On the other hand,  it has been shown,
using algebraic techniques,
in \cite{HM}
that, when $\cO$ is a Suzuki--Tits ovoid, $\cK$
cannot be obtained from a closed set of conics. In fact, in this
case the arc is not conical at all.
 
In order to
provide a direct representation of a Thas \cite{T74}
maximal arc in $\PG(2,q^2)$, where $q>4$ is an even prime power,
we 
shall use for $\PG(4,q)$ homogeneous coordinates
$(z,x_1,x_2,y_1,y_2)$.
The
hyperplane at infinity $\Sigma$ has equation $z=0$.
Let $\cS$ be a
regular spread of $\Sigma$ and denote by $\pi=\PG(2,q^2)$ the
corresponding Desarguesian plane obtained via Bruck--Bose
construction.
We shall use homogeneous coordinates $(z,x,y)$ for
$\pi$, in such a way that the line at infinity has equation $z=0$.
It is always possible to assume that, up to a projectivity,
the spread $\cS$ contains
the lines
\begin{equation}
\label{eq:ll}
\makeatletter
\begin{array}{r@{\langle}c@{,}c@{\rangle}}
 \ell_1= & (1,0,0,1) & (0,1,1,0) \\
 \ell_2= & (1,0,0,0) & (0,1,0,0) \\
 \ell_3= & (0,0,1,0) & (0,0,0,1) \\
 \end{array}.
\makeatother
\end{equation}
The map $\theta$ which realises the correspondence between the
points of $\PG(4,q)$ and those of $\PG(2,q^2)$ should map any line
$\ell$ of the spread $\cS$ in a point of $\PG(2,q^2)$. In
particular, to have
\[ \begin{array}{l@{=}l}
 \theta(\ell_1)&(0,1,1) \\
 \theta(\ell_2)&(0,1,0) \\
 \theta(\ell_3)&(0,0,1) \\
 \end{array} \]
we should choose
\[
 \theta:\left\{\begin{array}{l}
   \PG(4,q) \mapsto \PG(2,q^2) \\
   (z,x_1,x_2,y_1,y_2) \mapsto (z,x_1+\varepsilon x_2,\varepsilon y_1+y_2)
   \end{array}\right.,
\]
where $\varepsilon\in\GF(q^2)\setminus\GF(q)$ is a suitable element.

\fvset{numbers=left,frame=single,firstnumber=last}


\section{The code}
\label{sect:3}
In this section we
describe a GAP \cite{GAP4} program to construct a
Thas \cite{T74} maximal arc $\cK$ and determine a minimum degree
curve $\Gamma$ passing through all the points of $\cK$.
In our code it shall be constantly
assumed that  $q=2^{2t+1}$, with $t>1$
is a global variable.
\par
The simplest way to implement the geometry
$\PG(3,q)$ is to consider the
point orbit of $\GL(4,q)$ in its action on
left--normalised $4$--vectors.
\begin{Verbatim}
PG3:=Orbit(GL(4,q),[1,0,0,0]*Z(q)^0,OnLines);
\end{Verbatim}
\begin{remark}
It is often convenient
to represent the points of $\PG(3,q)$ as
integers in the range $1\ldots q^3+q^2+q+1$.
The number corresponding to any given point
is just the position of the corresponding normalised
vector in the list ${\tt PG3}$.
This is most interesting when ${\tt PG3}$ is generated
as the orbit of a point, say $(1,0,0,0)$,
under the action of a Singer group $\Theta$ of $\PG(3,q)$.
\end{remark}
We now introduce some utility functions.
 \begin{enumerate}
  \item {\tt LineAB} to compute the (projective) line
    over $\GF(q)$ through two points;
  \item {\tt LineAB2} to compute the (projective) line
    over $\GF(q^2)$ through two points.
  \item {\tt Conj} to get the conjugate of a point in
    $\PG(n,q^2)$ under the Frobenius morphism
    \[ x\mapsto x^q. \]
  \end{enumerate}
  \begin{Verbatim}
#Line (over GF(q))
LineAB:=function(a,b)
 return Set(Union([a],Set(GF(q),x->NormedRowVector(x*a+b))));
end;;

#Line (over GF(q^2))
LineAB2:=function(a,b)
 return Set(Union([a],Set(GF(q^2),x->NormedRowVector(x*a+b))));
end;;

#Conjugate of a point
Conj:=function(x)
 return(
 List(x,t->t^q));
end;;
\end{Verbatim}

The Suzuki group $\Sz(q)$ has two point orbits
in $\PG(3,q)$, of size respectively $q^3+q$ and
$q^2+1$. The latter is a Suzuki--Tits ovoid,
say $\cO={\tt Ov}$.
\begin{Verbatim}
Sg:=SuzukiGroup(IsMatrixGroup,q);
Or:=Orbits(Sg,PG3,OnLines);
Ov:=Filtered(Or,x->Size(x)=q^2+1)[1];
Ovp:=Set(Ov,x->Position(PG3,x));
\end{Verbatim}
We wrote ${\tt Ovp}$ for the set of all
points of ${\tt Ov}$
in the permutation representation.

The following code is used to write
the set $\Lambda={\tt AllLines}$
of all the lines of $\PG(3,q)$.
Since the full projective general linear
group $\PGL(4,q)$ is transitive on this set,
we may just consider the orbit of
\[ \ell_0=\left< (0,1,0,0), (1,0,0,0)\right> \]
under its action.

The group $\PGL(4,q)$ has to be
written as the action ${\tt Pgrp}$ of $GL(4,q)$
on normalised vectors. The line orbit is obtained
considering the action of this group ${\tt Pgrp}$
on the set of points, in the permutation
representation, of a given line.
\begin{Verbatim}
  L1:=LineAB([1,0,0,0]*Z(q)^0,[0,1,0,0]*Z(q)^0);
  L1p:=Set(L1,x->Position(PG3,x));
  Pgrp:=Action(GL(4,q),PG3,OnLines);
  AllLines:=Orbit(Pgrp,L1p,OnSets);
\end{Verbatim}


\begin{remark}
There might be more efficient ways to obtain the set $\Lambda$ as
union of line--orbits under the action of a Singer cycle $\Theta$ of
$\PG(3,q)$. In fact, see \cite{GLS}, the number to these
line--orbits is exactly  $q+1$ and a starter set for these (that is
a set of representatives for each of them) is given by all the lines
passing through a fixed point $P\not\in\cO^{+}$ tangent to the
elliptic quadric $\cO^+$ stabilised by the subgroup of order $q^2+1$
of $\Theta$.
\end{remark}

We are now in position to write the
set $T\cO={\tt TangentComplex}$
of all lines tangent to the ovoid ${\tt Ov}$.
This is simply done by enumerating the lines of
$\PG(3,q)$ which meet $\cO$ in just $1$ point.

The function ${\tt TCpx}$ is used to partition
this set according to the tangency point of the
lines themselves to the ovoid.
\begin{Verbatim}
TangentComplex:=
   Set(Filtered(AllLines,
               x->Size(Intersection(Ovp,x))=1),
       x->Set(x));

TCpx:=function(TC,O)
 return List(O,x->Filtered(TC,v->x in v));
end;;
\end{Verbatim}

As seen in Section \ref{sect:2},
given
three mutually skew lines $\ell_1,\ell_2,\ell_3$,
it is easy to write the opposite regulus $\cR^o$
they induce.
The regulus $\cR$ containing $L$ is then
obtained as $(\cR^o)^o$.

\begin{Verbatim}
#Functions to build up a
# regulus
# Here we use a permutation
# representation
OpRegulus:=function(a,b,c)
 return Filtered(AllLines,x->not(
   IsEmpty(Intersection(a,x)) or
   IsEmpty(Intersection(b,x)) or
   IsEmpty(Intersection(c,x))));
end;;

Regulus:=function(a,b,c)
 local l;
 l:=OpRegulus(a,b,c);
 return OpRegulus(l[1],l[2],l[3]);
end;;

#This function uses a normalised # vector representation
RegLines:=function(L)
 local Lp,Rp;
 Lp:=Set(L,x->Set(x,t->Position(PG3,t)));
 Rp:=Regulus(Lp[1],Lp[2],Lp[3]);
 return Set(Rp,
            x->Set(x,t->PG3[t]));
end;;
\end{Verbatim}

To construct a regular spread $\cS$ we
use the following functions:
\begin{enumerate}
\item {\tt LookForSpread0} which, given $4$ lines
 $\ell_1,\ell_2,\ell_3,\ell_4$, builds the set $R$ of
 all lines in reguli of the form $\cR(\ell_1,x,\ell_4)$ where
 $x\in\cR(\ell_1,\ell_2,\ell_3)\setminus\{\ell_1\}$;
\item {\tt RClosure} which determines $q^2+1$
  lines in the \emph{regular closure} of
  a set of lines $R$;
\item {\tt LookForSpread1}, {\tt LookForSpread2} and
  {\tt LookForSpread} which build the requested
  regular spread of  tangent lines to an ovoid ${\tt O}$.
\end{enumerate}

\begin{Verbatim}
# L = Set of 4 lines
LookForSpread0:=function(L)
 local Reg,RegT,x,Spr;
 Spr:=[];
 Reg:=Regulus(L[1],L[2],L[3]);
  for x in Difference(Reg,[L[1]]) do
   RegT:=Regulus(L[1],x,L[4]);
   Spr:=Union(Spr,RegT);
  od;
 return Spr;
end;;

RClosure0:=function(S)
 local x,X,R,V;
 X:=Combinations(S,3);
 R:=ShallowCopy(S);
 for x in X do
  R:=Union(R,Regulus(x[1],x[2],x[3]));
  if Size(R)=q^2+1 then return R;
  fi;
 od;
 return R;
end;;

RClosure:=function(S)
 local f,T;
 f:=false;
 T:=RClosure0(S);
 if not(T=S) then
  Print(Size(T),"-",Size(S),"\n");
  return RClosure(T);
 else
  Print("Closed\n");
  return T;
 fi;
end;;

# Hint for regulus
LookForSpread1:=function(TC,x,O)
local Tp,Ct,y,R1,S2,TC2;
  R1:=Regulus(x[1],x[2],x[3]);
  if not(IsSubset(TC,R1)) then return fail; fi;
  TC2:=Filtered(TC,x->IsEmpty(Intersection(x,Union(R1))));
  for y in TC2 do
   Print(".\n");
   S2:=LookForSpread0([x[1],x[2],x[3],y]);
   if IsSubset(TC,S2) then return (S2); fi;
  od;
  return fail;
end;;

LookForSpread2:=function(TC,O)
 local Tp,Ct,x,R;
 Tp:=Set(TCpx(TC,O),x->Set(x));
#First regulus
 Ct:=Filtered(Cartesian(Tp{[1..3]}),
              t->IsEmpty(Intersection(t[1],t[2])) and
                 IsEmpty(Intersection(t[1],t[3])) and
                 IsEmpty(Intersection(t[2],t[3])));
#Look for a second (compatible) regulus
 for x in Ct do
  R:=LookForSpread1(TC,x,O);
  if IsList(R) then return R; fi;
 od;
 return fail;
end;;

LookForSpread:=function(TC,O)
 local T;
 T:=RClosure(LookForSpread2(TC,O));
 if IsSubset(TC,T) then return T; fi;
 return fail;
end;;
\end{Verbatim}

To check if any given spread is regular, we
verify that it contains the regulus spanned by
any three of its elements.
\begin{Verbatim}
#Check if a spread is regular
IsRegularS:=function(S)
 local x,X,r;
 X:=Combinations(S,3);
 while Size(X)>2 do
 x:=X[1];
 r:=Regulus(x[1],x[2],x[3]);
  if not(IsSubset(S,r)) then
   Print(Size(Intersection(S,r)),"\n");
   return false;
  else
   X:=Difference(X,Combinations(r,3));
   Print(Size(X),"\n");
  fi;
 od;
 return true;
end;;
\end{Verbatim}

Our next step in constructing a model of $\PG(2,q^2)$
is to embed $\PG(3,q)$ in $\PG(4,q)$ as hyperplane
at infinity, as seen in Section \ref{sect:2}.
The function
${\tt EmbedPG3}$ does just this;
${\tt EmbedSpr}$ is a utility function
to embed sets of points of $\PG(3,q)$ in $\PG(4,q)$ and it
is most useful for spreads.
\begin{Verbatim}
# Embed PG(3,q) in PG(4,q) as
# hyperplane at infinity
EmbedPG3:=function(L)
 return Set(L,x->Concatenation([0*Z(q)],x));
end;;

EmbedSpr:=function(L)
 return Set(L,x->EmbedPG3(x));
end;;
\end{Verbatim}



Suppose now ${\tt Spr}$ to be a regular spread
of tangent lines to ${\tt Ov}$. We shall determine
a linear transformation $\mu$ of $\PG(3,q)$ such that
the spread $\mu({\tt Spr})$
contains the lines $\ell_1,\ell_2,\ell_3$
of \eqref{eq:ll}.
Recall that, for any spread $\cS$ of $\PG(3,q)$,
there exists a line $L_{\cS}$ of $\PG(3,q^2)\setminus\PG(3,q)$
such that
\[ \cS=\{ PP^q: P\in L_{\cS} \}. \]
Clearly, the spread $\cS$ is uniquely determined by the line
$L_{\cS}$, although different lines might be associated to the same
spread. The following function, ${\tt LookForLine}$, computes one of
these lines.
\begin{Verbatim}
LookForLine:=function(spr)
 local PSpr,xSpr,LLa,x, y, fl,xq;
 PSpr:=List(spr,x->LineAB2(PG3[x[1]],PG3[x[2]]));
 xSpr:=List(PSpr,x->Difference(x,PG3));
 LLa:=List(Cartesian(PSpr[1],PSpr[2]),x->LineAB2(x[1],x[2]));;
 for x in LLa do
  Print("x=",x[1],",",x[2],"\n");
#The lines should be disjoint from PG(3,q)
  if not(IsEmpty(Intersection(x,PG3))) then
   Print("!\n");
   continue;
  fi;
#They should also meet any component of the spread
  fl:=true;
  for y in xSpr do
   if IsEmpty(Intersection(y,x)) then
    Print("%");
    fl:=false;
    break;
   fi;
   Print(".");
  od;
  if not(fl) then continue; fi;
# The conjugate line
# should also meet any component of the spread
  xq:=Set(x,t->Conj(t));
  for y in xSpr do
   if IsEmpty(Intersection(y,xq)) then
    fl:=false;
    Print("^");
    break;
   fi;
  Print(",");
  od;
#If this is the case, then we have found
# what we were looking for
  if fl then return x; fi;
 od;
#Bad luck here.
return fail;
end;;
\end{Verbatim}
Denote now by ${\tt LCanon}$ the line of $\PG(3,q^2)$
associated with a spread, say ${\tt SCanon}$,
containing $\ell_1,\ell_2,\ell_3$.
\begin{Verbatim}
 GCanon:=[
 LineAB([1,0,0,1]*Z(q)^0,[0,1,1,0]*Z(q)^0),
 LineAB([1,0,0,0]*Z(q)^0,[0,1,0,0]*Z(q)^0),
 LineAB([0,0,1,0]*Z(q)^0,[0,0,0,1]*Z(q)^0)];
 GCanonP:=Set(GCanon,
     x->Set(x,t->Position(PG3,t)));
 RCanon:=Regulus(GCanonP[1],GCanonP[2],GCanonP[3]);
# look for a fourth line to generate the spread
 Get4th:=function(R)
  local j,L4;
  j:=1;
  repeat
   L4:=AllLines[j];
   j:=j+1;
   until IsEmpty(Intersection(L4,Union(R)));
  return L4;
 end;;
 L4:=Get4th(RCanon);
 SCanon:=RClosure(Union(GCanonP,[L4]));
 LCanon:=LookForLine(SCanon);
\end{Verbatim}
We may now actually determine a matrix of $\GL(4,q)$
inducing a collineation $\mu$ in $\PG(3,q)$
transforming the general spread ${\tt Spr}$ into ${\tt SCanon}$.
\begin{Verbatim}
 SprToCanon:=function(Spr)
  local Lx,M0,N0;
  Lx:=LookForLine(Spr);
  M0:=TransposedMat([Lx[1],Conj(Lx[1]),Lx[2],Conj(Lx[2])]);
  N0:=TransposedMat([LCanon[1],Conj(LCanon[1]),
                     LCanon[2],Conj(LCanon[2])]);
  return N0*M0^(-1);
 end;;
\end{Verbatim}
Let then ${\tt M}={\tt SprToCanon(Spr)}$ and suppose
${\tt SprT}=\mu({\tt Spr})$ and
${\tt OvT}=\mu({\tt Ov})$.
\begin{Verbatim}
# New spread
SprT:=Set(Spr,x->Set(x,t->NormedRowVector(M*t)));
#Consider also the image of the ovoid under the
# collineation induced by M
OvT:=Set(Ov,x->NormedRowVector(M*x));
\end{Verbatim}
We still need to determine the parameter $\varepsilon$ in the
correspondence $\theta:\PG(4,q)\mapsto\PG(2,q^2)$ of Section
\ref{sect:2}.
\begin{Verbatim}
PG4ToPG2:=function(P,eps)
 return NormedRowVector([P[1],P[2]+eps*P[3],eps*P[4]+P[5]]);
end;;
\end{Verbatim}
We
may proceed as follows.
\begin{Verbatim}
LookForEps:=function(Spr)
 local t,r,sp1,L,R1;
 L:=
 [LineAB([1,0,0,1],[0,1,1,0]),
  LineAB([1,0,0,0],[0,1,0,0]),
  LineAB([0,0,1,0],[0,0,0,1])]*Z(q)^0;
 R1:=RegLines(L);
 sp1:=Difference(Spr,R1);
 t:=sp1[1];
 r:=Filtered(
 Difference(Elements(GF(q^2)),Elements(GF(q))),
   eps->
   (t[1][1]+t[1][2]*eps)/(t[1][3]*eps+t[1][4])=
   (t[2][1]+t[2][2]*eps)/(t[2][3]*eps+t[2][4]));
 return r;
end;;

eps:=LookForEps(SprT)[1];
\end{Verbatim}

We are now in position to use \cite{T74} construction in
order to obtain a maximal arc.
We first embed $\PG(3,q)$ in $\PG(4,q)$ as the
hyperplane at infinity;
${\tt EOvT}$ is the image under this embedding of
the transformed ovoid (under the collineation
given by $\mu$);  then, we compute the \emph{affine}
cone ${\tt FullCone2}$ with vertex
\[ {\tt Vtx}=(1,0,0,0,0) \]
and basis ${\tt EOvT}$.
The image of this cone under $\theta={\tt PG4ToPG2}$ is
a maximal arc ${\tt Arc}$ of $\PG(2,q^2)$.
\begin{Verbatim}
# Embed OvT\subseteq PG(3,q) in PG(4,q)
EOvT:=EmbedPG3(OvT);
# ... and build the full cone in AG(4,q)
# with vertex
Vtx:=[1,0,0,0,0]*Z(q)^0;
# and basis OvT
FullCone:=Difference(Union(Set(EOvT,x->LineAB(x,Vtx))),EOvT);
# The requested maximal arc is the image of
# the cone under the map PG4ToPG2
Arc:=Set(FullCone,x->PG4ToPG2(x,eps));
\end{Verbatim}

We might want to check that we actually obtained an arc of degree
$q$. The function ${\tt CheckSecants}$,  whose parameter is a set ${\tt
X}$, verifies  that all of the secants of ${\tt X}$ meet ${\tt X}$
in $q$ points.
The function ${\tt CheckArc}$   performs the full check that:
 \begin{enumerate}
 \item all secants meet the given set in $q$ points;
 \item there is no tangent line through any point
  of ${\tt X}$.
 \end{enumerate}
\begin{Verbatim}
# Check if a set X is an arc
# step 0:
#  verify if all secants meet X in
#  q points
 CheckSecants0:=function(X)
 local C,l,XX;
 C:=Combinations(X,2);
 XX:=[];
 while(not(IsEmpty(C))) do
  l:=LineAB2(C[1][1],C[1][2]);
  if not(Size(Intersection(l,X))=q) then
   Print(Size(Intersection(l,X)),"\n");
   return [false,[]];
  fi;
  C:=Difference(C,Combinations(Intersection(l,X),2));
  Print("!",Size(C),"!\n");
  Add(XX,l);
 od;
 return [true,XX];
end;;

CheckSecants:=function(X)
 return (CheckSecants0(X)[1]);
end;;

CheckArc:=function(X)
 local C,l,XX,x;
 C:=Combinations(X,2);
#Computes all the secants;
 XX:=CheckSecants0(X);
 if not(XX[1]) then return false; fi;
 for x in X do
  l:=Filtered(XX[2],t->x in t);
  if Size(l)<q^2+1 then return false; fi;
 od;
 return true;
end;;
\end{Verbatim}

We are now ready to obtain a minimum degree curve covering
the arc $\cK={\tt Arc}$.
We proceed as follows:
\begin{enumerate}
 \item We determine all monic monomials
  in two variables of degree at most $i$ over $\GF(q^2)$.
  This is done by the function ${\tt AllMon}$.
\item A polynomial
  \[ f(x,y)=\sum_{i,j}c_{ij}x^iy^j \]
  describes a curve covering $\cA$ if, and only if,
  the coefficients $c_{ij}$ are a solution of the
  homogeneous linear system given by
  \begin{equation}
    \sum_{i,j}c_{ij}p_x^ip_y^j=0,\qquad P=(1,p_x,p_y)\in\cA;
  \end{equation}
\item Hence, we introduce a function,
  ${\tt BuildMat}$ that, for
  any given list of points $\cK$ and
  a maximum degree ${\tt i}$ generates
  the matrix whose rows are exactly the
  evaluations of the monomials in
  ${\tt AllMon(i)}$, computed on the second and
  third coordinate of any point in
  $\cK$.
  In other words, if
  \[ {\tt AllMon(i)}=\{f_1(x,y),f_2(x,y),\ldots,f_k(x,y)\} \]
  and $P=(1,p_x,p_y)\in{\cK}$, then the row of ${\tt BuildMat(K,i)}$
  corresponding to $P$ would be
  \[ [ f_1(p_x,p_y), f_2(p_x,p_y), \ldots, f_k(p_x,p_y) ]. \]
 \item If ${\tt BuildMat(K,i)}$ has full rank, then
  the only polynomial of degree at most ${\tt i}$
  in $x,y$ giving a curve
  covering all points of $\cK$ is
  the zero--polynomial.
\end{enumerate}

\begin{Verbatim}
RR:=PolynomialRing(GF(q^2),["x","y"]);
AllMon:=function(i)
 local l;
 l:=Filtered(Cartesian([0..i],[0..i]),t->t[1]+t[2]<i+1);
 return List(l,t->RR.1^t[1]*RR.2^t[2]);
end;;

BuildMat:=function(K,i)
 local m;
 m:=AllMon(i);
 return List(K,x->
        List(m,t->Value(t,[RR.1,RR.2],[x[2],x[3]])));
end;;
\end{Verbatim}
We need to know the
minimum index ${\tt i}$ such that
${\tt Buildmat(Arc,i)}$ has not
full rank.
The following function takes as
parameters the arc ${\tt K}$ and
a maximum degree to test.
Observe that
\[ \xi({\tt i})=\rank ({\tt BuildMat(K,i)})-|{\tt AllMon(i)}| \]
is non--increasing in ${\tt i}$.
Hence, to look for $i$, we may use an iterative approach:
consider an initial interval to test $[a\ldots b]$, let
$c=\lfloor\frac{a+b}{2}\rfloor$ and compute $\xi(c)$.
If $\xi(c)=0$,
then the first value ${\tt i}$ such that $\xi({\tt i})<0$
may possibly be found in $[(c+1)\ldots b]$; on the other hand, if
$\xi(c)<0$, such ${\tt i}$ is to be found in $[a\ldots c]$.
We keep bisecting the interval till it
contains just one value $c'$. If $\xi(c')<0$, then
${\tt i}=c'$ is returned; otherwise the algorithm fails.
\begin{Verbatim}
GetIndex:=function(A,mi)
 local tidx,c,d,r;
 tidx:=[1..mi];
 while(Size(tidx)>1) do
  c:=Int((tidx[1]+tidx[Size(tidx)])/2);
  d:=BuildMat(A,c);
  r:=Rank(d);
  Print("c=",c," t=",tidx,"\n");
  Print("r=",r," s=",Size(d[1]),"\n");
  if r=Size(d[1]) then
   tidx:=[(c+1)..tidx[Size(tidx)]];
  else
   tidx:=[tidx[1]..c];
  fi;
 od;
 Print(tidx,"\n");
 c:=tidx[1];
 d:=BuildMat(A,c);
 r:=Rank(d);
 if not(r=Size(d[1])) then
   return c;
 else
   return fail;
 fi;
end;;
\end{Verbatim}
Remark that the affine curve of equation
\[ (x^{q^2}-x)=0 \]
has degree $q^2$ and passes through all the points
of the affine plane $AG(2,q^2)$ (hence, also through all those of $\cK$).
Thus, this value may be chosen as the maximum degree $i$ to test
in ${\tt GetIndex}$.

\begin{Verbatim}
i:=GetIndex(Arc,q^2);
\end{Verbatim}

We may now obtain the coefficients of the polynomial giving the
curve  we are looking for, by just solving a linear system of
equations.

\begin{Verbatim}
MatOk:=BuildMat(Arc,i);;
SolV:=NullspaceMat(TransposedMat(MatOk))[1];
\end{Verbatim}

The values in ${\tt SolV}$ are now used to
write the equation of the
curve. This is done by the
function ${\tt VecToPoly}$.
\begin{Verbatim}
VecToPoly:=function(v,i)
 local m;
 m:=AllMon(i);
 return Sum(List([1..Size(v)],x->m[x]*v[x]));
end;;

pp:=VecToPoly(SolV,i);
\end{Verbatim}

\begin{remark}
When $q=8$, the construction of \cite{T74} gives
a $(456,8)$--maximal arc $\cK$ of $\PG(2,64)$.
If the ovoid $\cO$ chosen for this construction is
an elliptic quadric, then the minimum degree of
a curve $\Gamma$ containing all the points of $\cK$ is
$7$ and this curve splits into $3$ conics and
a line. On the other hand, if the Suzuki--Tits ovoid
is chosen, then the minimum degree of such a curve
$\Gamma$ is $22$ and it splits into an irreducible
curve of degree $17$, and $5$ lines.
\end{remark}

\bigskip
\begin{minipage}{6.8cm}
\begin{obeylines}
{\sc Angela Aguglia}
Dipartimento di Matematica
Politecnico di Bari
Via G. Amendola 126/B
70126 Bari
Italy
{\tt a.aguglia@poliba.it}
\end{obeylines}
\end{minipage}
\qquad
\begin{minipage}{6.8cm}
\begin{obeylines}
{\sc Luca Giuzzi}
Dipartimento di Matematica
Politecnico di Bari
Via G. Amendola 126/B
70126 Bari
Italy
{\tt l.giuzzi@poliba.it}
\end{obeylines}
\end{minipage}

\end{document}